\documentclass[12pt]{amsart}

\usepackage{amsmath}
\usepackage{amssymb}
\usepackage{amsthm}
\usepackage{graphicx} 
\usepackage{url}
\usepackage{hyperref}
\usepackage{multirow}
\usepackage{array}
\usepackage{color}

\theoremstyle{plain}
\newtheorem{theorem}{Theorem}

\theoremstyle{definition}

\newtheorem{remark}[theorem]{Remark}

\newcommand{\cir}{\mathrm{s}}
\newcommand{\cro}{\mathrm{c}}

\newcommand{\fr}{\mathrm{fr}}

\renewcommand{\ge}{\geqslant}

\newcommand{\inter}{\mathrm{i}}

\newcommand{\lk}{\mathrm{lk}}
\renewcommand{\le}{\leqslant}

\newcommand{\NN}{\mathbb{N}}

\newcommand{\RR}{\mathbb{R}}

\newcommand{\SeifS}{\mathrm{G}_S}
\newcommand{\SeifSD}{\mathrm{G}_{S_D}}
\newcommand{\Sph}{\mathbb{S}}

\newcommand{\ZZ}{\mathbb{Z}}

\title{Signature and concordance of positive knots}

\author{S.~Baader}
\address{Mathematisches Institut, Universit\"at Bern, Sidlerstrasse 5, 3012 Bern, Schweiz}
\email{sebastian.baader@math.unibe.ch}
\author{P.~Dehornoy}
\address{Univ. Grenoble Alpes, CNRS, IF, F-38000 Grenoble, France}
\email{pierre.dehornoy@ujf-grenoble-alpes.fr}
\urladdr{\url{http://www-fourier.ujf-grenoble.fr/~dehornop/}}
\author{L.~Liechti}
\address{Institut de Math\'ematiques de Jussieu, Universit\'e Pierre et Marie Curie, 75252 Paris Cedex 05, France}
\email{livio.liechti@imj-prg.fr}
\thanks{The third author is supported by the Swiss National Science Foundation (project no.\ 137548).}

\begin{document}
\begin{abstract}
We derive a linear estimate of the signature of positive knots, in terms of their genus. As an application, we show that every knot concordance class contains at most finitely many positive knots.
\end{abstract}
\maketitle
\section{Introduction}
 Algebraic links arise in the context of plane curve singularities. In terms of knot theory, they are certain iterated cables of torus links. 
 Their classification is well-understood, even on the level of smooth concordance: concordant algebraic links are equal~\cite{Li}. 
 This feature is believed to extend to various larger classes of knots, in particular to closures of positive braids. 
 As pointed out by Baker, this would follow from the Slice-Ribbon Conjecture~\cite{Ba}. 
 Stoimenow proposed a weaker formulation that is supposed to hold for the much larger class of positive knots and proved it for closures of positive braids~\cite{Stoimenow}, 
 as well as special alternating knots~\cite{Stoimenow2}. As we will see, an extension of his technique works for positive knots.

\begin{theorem}\label{T1}
Every topological, locally-flat knot concordance class contains at most finitely many positive knots.
\end{theorem}

No larger class of knots is likely to share this feature. Rudolph's construction of transverse $\mathbb{C}$-links yields infinite families of smoothly slice quasipositive knots~\cite{Rud2}. Furthermore, Baker describes infinite families of smoothly concordant (in fact, ribbon concordant) strongly quasipositive knots~\cite{Ba}, based on Hedden's work on Whitehead doubles~\cite{He}. Our proof relies on a careful analysis of the signature function of positive knots. The main ingredient is a linear estimate of the signature invariant $\sigma$ of positive links.

\begin{theorem}\label{T}
The following inequalities hold for all positive links $L$: 
\[ \frac{1}{24} b_1(L) \le \sigma(L) \le b_1(L).\]
\end{theorem}

Here $b_1(L)$ denotes the minimal first Betti number of all Seifert surfaces of the link~$L$. The earliest result in this direction is due to Rudolph~\cite{Rudolph}: 
closures of positive braids have positive signature\footnote{There exist two opposite conventions concerning the sign of the signature. Here we adopt Rudolph's convention: 
positive links have positive signature.}. This was independently extended to positive knots by Cochran and Gompf~\cite{CG} and to positive links by Traczyk and Przytycki~\cite{Traczyk, Prz}. 
Later Stoimenow improved Rudolph's result by showing that the signature is bounded by an increasing function of the first Betti number~\cite{Stoimenow}. 
The first linear bound for positive braid links was derived by Feller~\cite{Peschito}: $$\frac{1}{100} b_1 \le \sigma \le b_1.$$
Recently, the lower bound has been improved to $\frac{1}{8}b_1$ by the same author~\cite{Pegascito}. He conjectures the optimal bound to be $\frac{1}{2}b_1$ and proves this for closures of positive $4$-braids. 
An extension of this conjecture to the larger class of positive links is conceivable, as confirmed by examples with small crossing number.

Theorem~2 yields a linear bound for the topological 4-genus $g_4$ of positive knots, in terms of their genus: \[ \frac{1}{24} g(K) \le g_4(K). \]
This follows from Kauffman and Taylor's signature bound~\cite{KT}, $\sigma \le 2g_4(K)$, combined with the equality $2g(K)=b_1(K)$, valid for all knots $K$. Again, this does not extend to (strongly) quasipositive knots, thanks to a result by Rudolph~\cite{Rud2}: every Seifert form can be realized by a quasipositive surface. In particular, there exist topologically slice, strongly quasipositive knots of arbitrarily high genus.

The following section contains basic facts about the signature function, non-orientable spanning surfaces and the signature formula of Gordon and Litherland~\cite{GL}. The proofs of Theorems~1 and~2 are contained in Sections~4 and~3, respectively. 

\medskip
{\sc Acknowledgements.}
We would like to thank Peter Feller for all the inspiring discussions we have had.


\section{The signature function}

Let $L$ be a link and $\omega \in \Sph^1$. 
Levine and Tristram defined the $\omega$-{\it signature} $\sigma_\omega(L)$ of $L$ to be the signature of the hermitian matrix $$M_\omega=(1-\omega)A + (1-\bar\omega)A^t,$$
where $A$ is a Seifert matrix for the link $L$~\cite{Levine, Tristram}.
This does not depend on the choice of Seifert matrix $A$ and for $\omega = -1$, it equals Trotter's definition of the classical signature invariant $\sigma(L)$~\cite{Trotter}.
Furthermore, since $$\text{det}(M_\omega) = \text{det}(-(1-\bar\omega)(\omega A - A^t)) = \text{det}(-(1-\bar\omega))\Delta_L(\omega),$$
the signature function is locally constant except at zeroes of the Alexander polynomial $\Delta_L$, i.e.\ at finitely many points.
A crossing of an oriented link diagram is \emph{positive} if, when following the bottom strand, the top strand comes from the left, and \emph{negative} otherwise (see Fig.~\ref{F:Type}).
An oriented link is \emph{positive} if it admits a diagram with only positive crossings.
We will use that changing a positive crossing to a negative one does not increase the $\omega$-signature for any $\omega \in \Sph^1$. 

\subsection{Gordon and Litherland's formula for the signature $\sigma(L)$}
\label{signature_section}
Given a surface~$S$ (orientable or not) embedded in the 3-sphere and a curve~$x$ on~$S$, the \emph{normal displacement} $x^\pm$ of~$x$ is the multi-curve in~$\Sph^3\setminus S$ obtained by slightly pushing $x$ normally off the surface in both directions. 
It is connected if and only if $x$ is connected and has a non-orientable annular neighbourhood. 
The \emph{Goeritz form} of~$S$ is the symmetric bilinear form~$\SeifS$ on~$H_1(S, \ZZ)$ defined by~$\SeifS(x,y)=\lk(x^\pm,y)$~\cite{Goeritz}. 
If $S$ is orientable, the Goeritz form coincides with the symmetrised Seifert form.

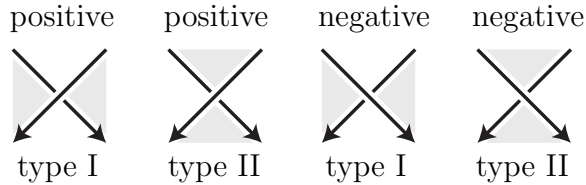
\begin{figure}[htb]
	\begin{picture}(70,25)(0,0)
		\put(0,4){\includegraphics[width=7.5cm]{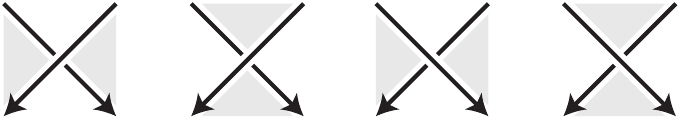}}
		\put(0,20){positive}
		\put(20.5,20){positive}
		\put(41,20){negative}
		\put(61.5,20){negative}
		\put(1,0){type I}
		\put(21,0){type II}
		\put(42,0){type I}
		\put(62,0){type II}
	\end{picture}
	\caption{\small The sign of a crossing depends on the orientation of the strands and the type depends on whether the orientations of the two strands both agree with a local orientation of the spanning surface.}
	\label{F:Type}
\end{figure}

\begin{figure}[htb]
	\includegraphics[width=.9\textwidth]{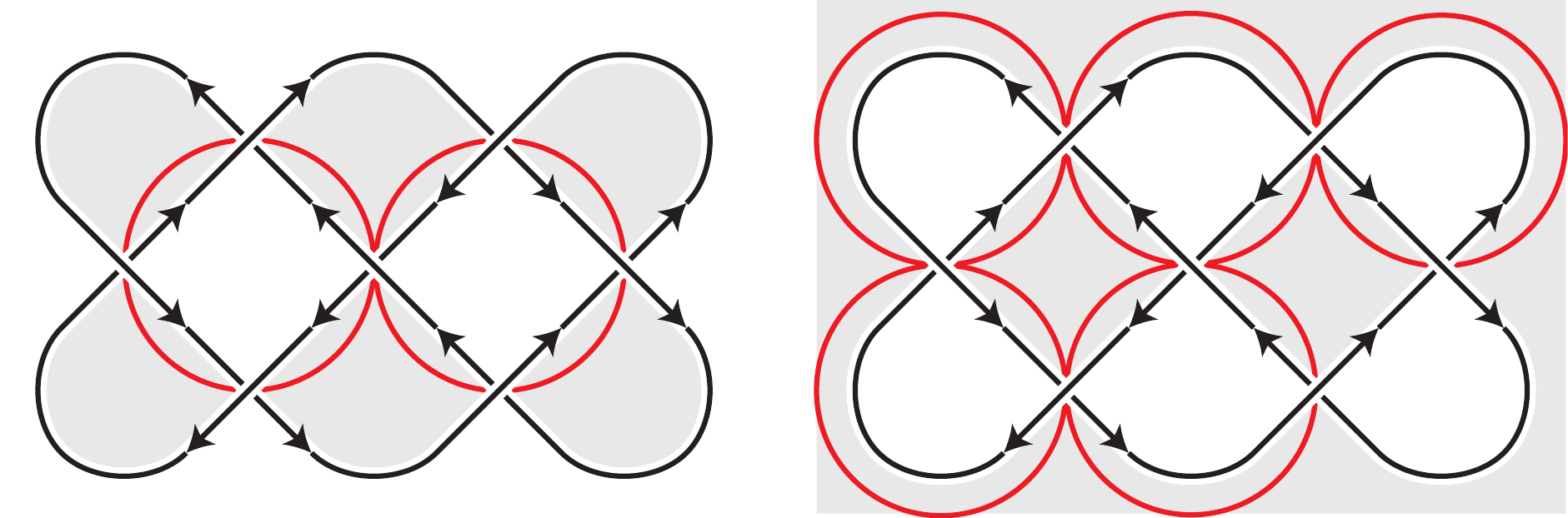}
	\caption{\small The two checkerboard surfaces associated to a positive planar link diagram and two bases of their first homology group (in red).}
	\label{F:BW}
\end{figure}

The complement of every link diagram in the plane admits two checkerboard colourings. 
To each of them is associated a spanning surface for the link obtained by lifting the black faces to 3-space (see Fig.~\ref{F:BW}). 
Since they need not be orientable, the surfaces may have much smaller genus than the link. 
Building on the work of Goeritz~\cite{Goeritz}, Gordon and Litherland gave a formula for computing the signature of the link using this data only~\cite{GL}. 
For an oriented link diagram $D$ with a given checkerboard colouring of~$\RR^2\setminus D$, we denote by~$S_D$ the associated checkerboard surface. 
A crossing of $D$ is of \emph{type~I} if the orientation of the two strands induce the same local orientation on~$S_D$ and of \emph{type~II} otherwise (see Fig.~\ref{F:Type}). 
The formula of Gordon and Litherland, stated below, contains a quantity $\mu(S_D)$. In the case of positive diagrams $D$, it is simply the number of type~II crossings.

\begin{theorem}[compare~\cite{GL} Thm $6''$]
\label{gordonlitherland}
	For $D$ a positive link diagram and $S_D$ the surface associated to a checkerboard colouring, the signature of the associated link is equal to~$\mathrm{sign}(\SeifSD) + \mu(S_D)$.
\end{theorem}

In order to compute the signature of the Goeritz form, we assume that the diagram $D$ is \emph{reduced}.
This means that smoothing any crossing of $D$ does not split off a link component.
Furthermore, we use a particular basis of~$H_1(S, \ZZ)$: 
for every white face of the checkerboard colouring except one, there is a generator $\gamma$ that runs around it (see Fig.~\ref{F:BW}).
For such a generator~$\gamma$, we denote by $\fr(\gamma)$ the number of type~I crossings minus the number of type~II crossings along~$\gamma$. 
For two generators $\gamma$ and $\gamma'$, we denote by~$\inter(\gamma, \gamma')$ their signed number of intersection, that is, the number of type I common crossing points minus the number of type II common crossing points.
Reducedness guarantees that no generator $\gamma$ runs twice through the same crossing.
In this basis, the coefficients of the Goeritz form are given by $(\SeifSD)_{i,i} = \fr(\gamma_i)$ for $i=j$ and $(\SeifSD)_{i,j} = -\inter(\gamma_i, \gamma_j)$ otherwise~\cite{Goeritz}.

\subsection{The first Betti number of a positive link}
Given a positive, non-split diagram~$D$ of a link $L$, Seifert's algorithm yields a canonical orientable surface~$S_{D}$. 
Cromwell proved that this surface is genus-minimising~\cite{Cromwell}. 
Let~$\cir(D)$ be the number of Seifert circles and let $\cro(D)$ be the number of crossings of~$D$. By Cromwell's theorem, we have $b_1(L) =  \cro(D)-\cir(D)+1.$

\section{Proof of Theorem~\ref{T}}

\begin{figure}[htb]
	\includegraphics[width=.6\textwidth]{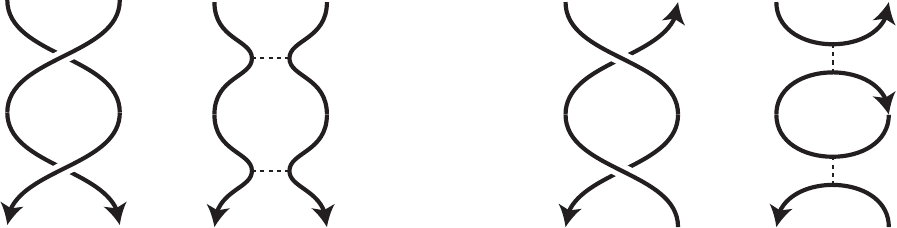}
	\caption{\small A face with two edges may (right) or may not (left) produce a Seifert circle.}
	\label{F:Circle}
\end{figure}

\begin{figure}[htb]
	\includegraphics[width=.75\textwidth]{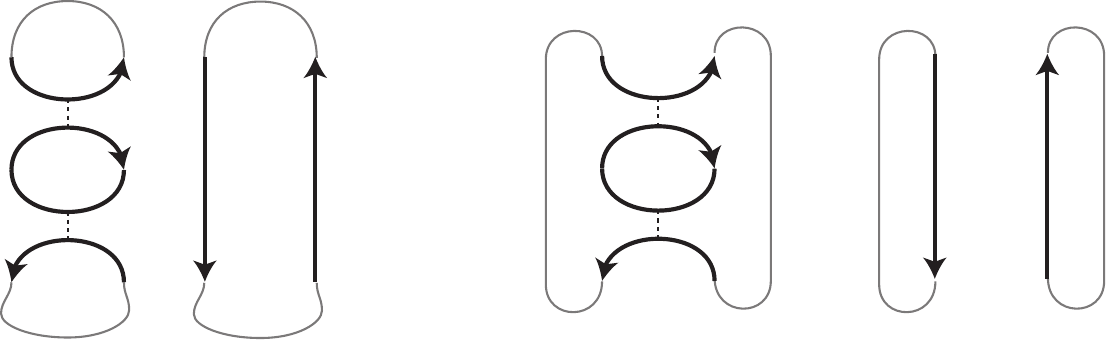}
	\caption{\small Untwisting the full twist given by two crossings that would yield a Seifert circle either reduces by two the number of crossings and the number of Seifert circles (left) or only reduces the number of crossings (right).}
	\label{F:Reducing}
\end{figure}

Let $D$ be a reduced positive diagram of a link $L$. 
In particular, there are no faces with only one edge.
We will modify the diagram $D$ at certain faces with two edges. 
There are two types of these depending on whether the boundary builds a Seifert circle or not  (see Fig.~\ref{F:Circle}).
If yes, there are still two distinct possibilities. 
Untwisting the full twist given by both crossings does or does not reduce the first Betti number (see Fig.~\ref{F:Reducing}).
Actually, by the formula $b_1(L) = \cro(D)-\cir(D)+1$,
untwisting does not reduce the genus if and only if the Seifert circle in question is connected to two distinct other Seifert circles. 
We choose to untwist those full twists where untwisting does not reduce the first Betti number. 
By this procedure, we obtain a reduced diagram $D'=D(L')$ of a new positive link $L'$ such that $b_1(L)=b_1(L')$. 
Furthermore we have $\sigma(L) \ge \sigma(L')$, since the untwisting can be realised by positive-to-negative crossing changes.
We now show $\sigma(L') \ge \frac{1}{24} b_1(L')$ for the new link $L'$. This implies the same inequality for the link $L$.

\begin{figure}[htb]
	\begin{picture}(30,28)(0,0)
		\put(0,0){\includegraphics[width=2.7cm]{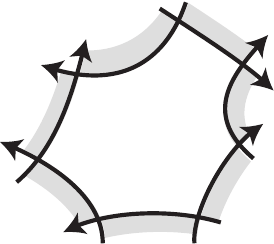}}
		\put(-2,6){II}
		\put(5,19){I}
		\put(16,24){I}
		\put(27,17){I}
		\put(26,8){II}
		\put(20,-2){I}
		\put(7,-3){II}
	\end{picture}
	\caption{\small Since type I crossings reverse the orientation around a white face, there is an even number of them.}
	\label{F:meven}
\end{figure}
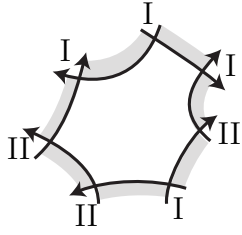

Choose a checkerboard colouring of the diagram $D'$ and let $S=S_{D'}$ be the surface defined by the black regions.
We further distinguish the basis curves $\{\gamma_i\}$ of $H_1(S,\ZZ)$ described above. 
We say that the curve $\gamma_i$ is of type $(m,n)$ if it runs through $m$ crossings of type I and $n$ crossings of type II, so $\fr(\gamma_i)=m-n$.
By considering the boundary orientation along a white face, we get that $m$ is always even (see Fig.~\ref{F:meven}). 
By our simplifications of the link diagram, we removed all curves of type $(0,1)$. 
Furthermore, all curves $\gamma_i$ of type $(0,2)$ correspond to Seifert circles with two edges we did not untwist before.
If two such Seifert circles meet at a crossing, then they meet at two crossings (otherwise untwisting a corresponding full twist would not reduce the first Betti number) and the corresponding part of~$D'$ corresponds to a 2-component Hopf link split from the rest of the link. 
Hence we can assume without loss of generality that no two such Seifert circles meet at a crossing, 
so that the number $\mu=\mu(S_{D'})$ of crossings of type II is at least twice the number of curves $\gamma_i$ of type $(0,2)$.

Writing $\gamma(m,n)$ for the number of curves $\gamma_i$ of type $(m,n)$, we get that $\gamma(0,1) = \gamma(1,n) = 0$ and $\gamma(0,2) \le \frac{1}{2} \mu.$
Furthermore, writing $\gamma_{<0}$ for the number of $\gamma_i$ with $\fr(\gamma_i)<0$, we have 
$$\gamma_{<0} = \gamma(0,2) + \sum_{\substack{n>m\\n>2}}\gamma(m,n).$$
Since every crossing of type II is met by at most two curves $\gamma_i$, we get 
$$2\mu \ge 2\gamma(0,2) + \sum_{\substack{n>m\\n>2}}n\gamma(m,n) \ge 2\gamma(0,2) + \sum_{\substack{n>m\\n>2}}3\gamma(m,n)$$
and, by combining with the equation above, we obtain
\begin{align*}
3\gamma_{<0} &=  2\gamma(0,2) + \sum_{\substack{n>m\\n>2}}3\gamma(m,n) + \gamma(0,2)\\
&\le 2\mu + \gamma(0,2)\\
&\le 2\mu + \frac{1}{2}\mu\\
&= \frac{5}{2}\mu,
\end{align*}
and therefore $\gamma_{<0} \le \frac{5}{6}\mu.$

Now we estimate the signature of $L'$ using the formula of Gordon and Litherland (Theorem~\ref{gordonlitherland}) discussed in Section~\ref{signature_section}. 
We actually need to estimate $\text{sign}(G_S)$.
In order to do this, we consider the curves $\gamma_i$ with $\fr(\gamma_i)\ge 0$. 
Let $\Gamma$ be the graph with a vertex $v_i$ for each curve $\gamma_i$ with $\fr(\gamma_i)\ge 0$ and an edge between two distinct vertices $v_i$ and $v_j$ if and only if $\inter(\gamma_i,\gamma_j)\ne0$.
Being defined via a planar construction, the graph $\Gamma$ is planar.
From the four-colour theorem\footnote{One could also use five-colourability instead of four-colourability. This would cause a drop in the constant of Theorem~\ref{T} to $\frac{1}{30}$.} 
it follows that one can choose a bipartite subgraph $\Gamma'\subset\Gamma$ containing at least one half of the vertices~\cite{AppelHaken}.
The Goeritz form restricted to the subspace generated by the curves $\gamma_i$ corresponding to the vertices of $\Gamma'$ is given by a block matrix $$\begin{pmatrix}
  D_1 & X\\
  X^{\top} & D_2
 \end{pmatrix},$$
where $D_1$ and $D_2$ are non-negative diagonal matrices. Since the signature of such a matrix is non-negative,
this gives a subspace of dimension at least $$\frac{1}{2}\gamma_{\ge0} = \frac{1}{2}((f_w -1) - \gamma_{<0})$$ restricted to which the Goeritz form $\SeifS$ has non-negative signature,
where $f_w$ is the number of white regions in the checkerboard colouring of $\RR^2\setminus D'$, so $\text{dim}(H_1(S,\ZZ))=f_w-1$.
Plugging this into the formula of Gordon and Litherland, we obtain
\begin{align*}
\sigma(L') &= \mu + \text{sign}(\SeifS)\\
&\ge \mu - (f_w -1) + \frac{1}{2}((f_w -1) - \gamma_{<0})\\
&\ge \mu - (f_w -1) + \frac{1}{2}(f_w -1) - \frac{1}{2}(\frac{5}{6}\mu)\\
&= \frac{7}{12}\mu - \frac{1}{2}(f_w -1).
\end{align*}

Using the other checkerboard surface, we have that white faces become black faces and vice versa, and crossings of type I become crossings of type II and vice versa. 
By exactly the same procedure, we get
$$\sigma(L') \ge \frac{7}{12}(\cro(D') - \mu) - \frac{1}{2}(f_b -1).$$
Finally, we sum both inequalities for $\sigma(L')$, observe $f_w+f_b-2=\cro(D')$, use $\cro(D') \ge b_1(L')$ and obtain 
\begin{align*}
2\sigma(L') &\ge \frac{7}{12}(\mu + \cro(D') - \mu) - \frac{1}{2}(f_w + f_b - 2)\\
&= \frac{7}{12} \cro(D') - \frac{1}{2} \cro(D')\\
&= \frac{1}{12}  \cro(D') \\
&\ge \frac{1}{12}  b_1(L').
\end{align*}

\section{Proof of theorem~\ref{T1}}

\begin{figure}[h]
  \def\svgwidth{320pt}
  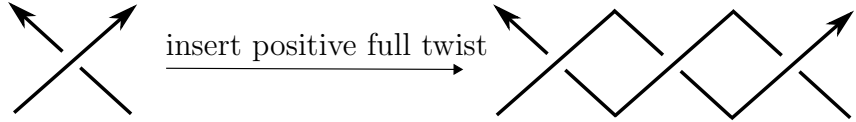
  \caption{Inserting a positive full twist at a crossing. Note that this operation can be reversed by a positive-to-negative crossing change.}
  \label{postwist}
\end{figure}

Suppose that there exists a topological, locally-flat concordance class $\mathcal{K}$ containing infinitely many positive knots $K_i$.
The average signature function depends only on the topological, locally-flat concordance class of a knot, so all the $K_i$ have identical signature function outside the zeroes of their Alexander polynomials~\cite{Po}.
In particular, all the $K_i$ have the same signature $\sigma = \sigma(\mathcal{K}) \in \NN$.
Thus, by Theorem~\ref{T}, the genera of the knots $K_i$ are bounded from above by $12\sigma$.
It then follows from a theorem of Stoimenow that there exist finitely many positive knot diagrams $D_j$, 
such that every knot $K_i$ is obtained from one of the $D_j$ by inserting a certain number of positive full twists at the crossings of $D_j$ (see Fig.~\ref{postwist})~\cite[Theorem 3.1]{Stoimenow01}.

\begin{figure}[h]
  \def\svgwidth{225pt}
  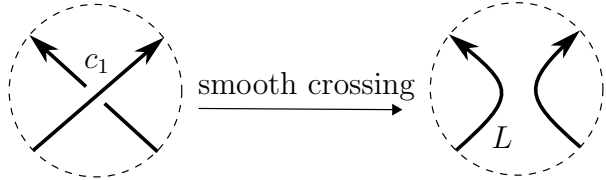
  \caption{The oriented smoothing of a crossing of a positive knot diagram yields a positive link diagram.}
  \label{smoothingc1}
\end{figure}

Since there are infinitely many $K_i$ but only finitely many $D_j$, 
we can assume without loss of generality that all the knots $K_i$ are obtained from one single reduced diagram $D$ by inserting a certain number of positive full twists at the crossings  $c_1, \dots , c_n$ of $D$.
Again, since there are infinitely many $K_i$ but only finitely many crossings, we can assume that the number of inserted full twists at one of the crossings, say $c_1$, becomes arbitrarily large as $i$ tends to infinity.
Now let $K(N)$ be the knot obtained by starting from the diagram $D$ and inserting $N$ positive full twists at $c_1$ and let $L$ be the link obtained from $D$ by smoothing the crossing $c_1$ as in Fig.~\ref{smoothingc1}.

Since the number of inserted positive full twists at the crossing $c_1$ becomes arbitrarily large as $i$ tends to infinity, for every $N \in \NN$, 
the knot $K(N)$ can be obtained from one of the $K_i$ by applying some positive-to-negative crossing changes. 
In particular, we have $$\sigma_{\omega}(\mathcal{K}) \ge \sigma_\omega(K(N))$$ for all non-algebraic $\omega \in \Sph^1$ and $N\in\NN$.
Here and in the following, we consider non-algebraic numbers $\omega \in \Sph^1$, guaranteeing that $\sigma_\omega$ as defined in Section~2 is a concordance invariant of knots.
Since $\mathcal{K}$ is the concordance class of a knot, $\sigma_{\omega}(\mathcal{K})=0$ for $\omega$ close enough to $1\in \Sph^1$. 
Choose $\omega_0$ such that $\sigma_{\omega_0}(\mathcal{K})=0$ and $\Delta_L(\omega_0) \ne 0$. 
The second condition can be achieved since the coefficient of the linear term of the Conway polynomial of a two-component link equals the linking number of its components 
and thus the Alexander polynomial $\Delta_L$ is not identically zero. 
Here we use that $L$ is non-split (otherwise $D$ would not be reduced) and positive. 
We now show that $N$ can be chosen so that $\sigma_{\omega_0}(K(N)) > 0$. This will be the contradiction we are looking for.

Since a Seifert surface for $K(N)$ can be obtained from a Seifert surface of $L$ by glueing a suitably twisted ribbon, a Seifert matrix for $K(N)$ can be chosen to be 
$$A = \begin{pmatrix}
x & w\\
v & B
\end{pmatrix},$$
where $B$ is a Seifert matrix for $L$, $v$ is a column vector, $w$ is a row vector and $x$ is a natural number. In fact, nothing depends on $N$ except for $x$ and $x > \frac{N}{2}$. 
We are interested in $\sigma_{\omega_0}(K(N))$, which is the signature of the matrix 
$$S_N = \begin{pmatrix}
(2-2\text{Re}(\omega_0))x & (1-\omega_0)w + (1 - \bar \omega_0)v^t\\
(1-\omega_0)v + (1 - \bar \omega_0)w^t & (1-\omega_0)A + (1 - \bar \omega_0)A^t
\end{pmatrix} =: 
\begin{pmatrix}
x' & \ast \\
\ast & B'
\end{pmatrix}.$$
The determinant of this matrix equals $x'\text{det}(B') + r$, where $r \in \NN$ does not depend on $N$ and $\text{det}(B') \ne 0$ since $\Delta_L(\omega_0) \ne 0$. 
From this it follows that for large enough $N$, the determinant of $S_N$ has the same sign as the determinant of $B'$. 
In particular, the matrix $S_N$ has one more positive eigenvalue than $B'$.
A theorem by Przytycki implies $\sigma_{\omega_0}(L) \ge 0$ since $L$ is a positive link~\cite{Prz}. Therefore the signature of $S_N$ is strictly positive for $N$ large enough.
In particular, this implies $\sigma_{\omega_0}(\mathcal{K}) > 0$, a contradiction.

\begin{remark}
 For positive knots, the smooth 4-genus equals the Seifert genus, by a theorem of Kronheimer and Mrowka~\cite{K-M}. 
 Using this, our proof of Theorem~\ref{T1} works in the smooth category without referring to Theorem~\ref{T}. 
\end{remark}


\bibliographystyle{siam}

\end{document}